\documentclass[12pt]{article}
\usepackage[russian]{babel}
\usepackage{amsthm,amsfonts,amssymb,amscd}

\begin{document}

%\begin{flushleft}
%ÓÄÊ 512.7
%\end{flushleft}

\begin{center}
\large \bf Factorial hypersurfaces
\end{center}\vspace{0.5cm}

\centerline{A.V.Pukhlikov}\vspace{0.5cm}

\parshape=1
3cm 10cm \noindent {\small \quad\quad\quad \quad\quad\quad\quad
\quad\quad\quad {\bf }\newline In this paper the codimension of
the complement to the set of factorial hypersurfaces of degree $d$
in ${\mathbb P}^N$ is estimated for $d\geqslant 4$, $N\geqslant
7$.

Bibliography: 12 titles.} \vspace{1cm}

\vspace{1cm}

\section*{Introduction}

{\bf 0.1. Statement of the main result.} A singular point $o$ of
an algebraic variety $V$ is {\it factorial}, if every prime
divisor $D\ni o$ in a neighborhood of this point is given by a
single equation $f\in{\cal O}_{o,V}$. (For a non-singular point
$o\in V$ this is a well known theorem of the classical algebraic
geometry.) Majority of the modern techniques work for factorial or
$\mathbb Q$-factorial varieties (in the latter case it is required
that some multiple of every prime divisor $D\ni o$ were given by
one equation, see, for instance, \cite{Me04} or any paper on the
minimal model program). A standard example of a non-${\mathbb
Q}$-factorial (and the more so, non-factorial) variety is the
three-dimensional cone in ${\mathbb P}^4$ over a non-singular
quadric in ${\mathbb P}^3$: no multiple of a plane passing through
the vertex of the cone (and contained in that cone) can be given
by one equation in the local ring of the vertex of the cone. The
aim of this paper is to estimate from below the codimension of the
complement to the set of factorial hypersurfaces of degree $d$ in
${\mathbb P}^N$, $N\geqslant 7$, $d\geqslant 4$. More precisely,
let ${\cal P}_{d,N+1}$ be the space of homogeneous polynomials of
degree $d$ in the coordinates $x_0,\dots,x_N$ on ${\mathbb P}^N$.
Let ${\cal P}^{\rm fact}_{d,N+1}\subset {\cal P}_{d,N+1}$ be the
subset, consisting of such $f\in {\cal P}_{d,N+1}$ that the
hypersurface $\{f=0\}$ is irreducible, reduced and
factorial.\vspace{0.1cm}

{\bf Theorem 0.1.} (i) {\it Assume that $4\leqslant d\leqslant N$
(that is, $\{f=0\}$ is a Fano hypersurface). Then the estimate
$$
\mathop{\rm codim}\left({\cal P}_{d,N+1}\setminus {\cal P}^{\rm
fact}_{d,N+1}\right)\geqslant \mathop{\rm min} \left[
3{d+N-5\choose N-2}-N,\,\, 5{d+N-6\choose N-3}\right]
$$
holds, and in the case $d=N$ (that is, $\{f=0\}$ is a Fano
hypersurface of index one) in the right hand side one can leave
only} $5{d+N-6\choose N-3}$.\vspace{0.1cm}

(ii) {\it Assume that $d\geqslant 2N$ (in particular, the
hypersurface $\{f=0\}$ is a variety of general type). Then the
following estimate holds:}
$$
\mathop{\rm codim}\left({\cal P}_{d,N+1}\setminus {\cal P}^{\rm
fact}_{d,N+1}\right)\geqslant {d+N-4\choose N-4}+4{d+N-5\choose
N-4}-4(N-3).
$$

In fact, we will obtain an estimate for the codimension of the
complement to the set ${\cal P}^{\rm fact}_{d,N+1}$ for any values
$d\geqslant 4$ (Theorem 3.1), just in the Fano case ($d\leqslant
N$) and in the case $d\geqslant 2N$ that estimate can be
essentially simplified to the inequalities of Theorem
0.1.\vspace{0.3cm}

%%%%%%%%%%%%%%%%%%%%%%%%%%%%%%%%%%%%%%%%%%%%%%%%%%%%%%%%%%%%%%%%
%%%%%%%%%%%% subsection 0.2

{\bf 0.2. The plan of the proof and the structure of the paper.}
Our proof is based on the famous Grothendieck's theorem \cite{CL}
and the technique of estimating the codimension of the set of
hypersurfaces in ${\mathbb P}^N$ with a singular set of prescribed
dimension, developed in \cite{Pukh15a,Pukh15b}. Grothendieck's
theorem claims that a variety with hypersurface singularities (in
fact, with complete intersection singularities) is factorial, if
the singular locus has codimension at least 4. For that reason, in
order to estimate the codimension of the complement to ${\cal
P}^{\rm fact}_{d,N+1}$, it is sufficient to estimate the
codimension of the subset consisting of polynomials $f\in {\cal
P}_{d,N+1}$ such that the singular locus of the hypersurface
$\{f=0\}$ has codimension at most 3. This is what we will do in
the present paper.\vspace{0.1cm}

The paper is organized in the following way. In \S 1 we compute
the codimensions of two sets of polynomials $f\in {\cal
P}_{d,N+1}$: such that the hypersurface $\{f=0\}$ has a linear
space of singular points and such that the hypersurface $\{f=0\}$
has a subvariety of singular points which is a hypersurface in a
linear space. After that we set the general problem of estimating
the codimension of the set of polynomials $f\in {\cal P}_{d,N+1}$,
such that the singular locus of the hypersurface $\{f=0\}$ is of
dimension at least $i\geqslant 1$, and state the main technical
result --- Theorem 1.1, solving this problem.\vspace{0.1cm}

\S 2 contains the proof of Theorem 1.1 by means of the technique
developed in \cite{Pukh15a,Pukh15b}. In \S 3 we obtain an estimate
of the codimension of the set of non-factorial hypersurfaces,
which implies Theorem 0.1.\vspace{0.3cm}

%%%%%%%%%%%%%%%%%%%%%%%%%%%%%%%%%%%%%%%%%%%%%%%%%%%%%%%%%%%%%%%%
%%%%%%%%%%%% subsection 0.3

{\bf 0.3. Historical remarks and acknowledgements.} Factoriality
of algebraic varieties is a very old topic in Algebraic Geometry,
with lots of papers written on the subject. We will only point out
a few recent papers that demonstrate that the topic is still
actively investigated today:
\cite{Ch06f,Ch10f,ChPark06,Chil04a,Polizzi14}. Various technical
points related to the constructions of this paper were discussed
by the author in his talks given in 2009-2014 at Steklov
Mathematical Institute. The author thanks the members of divisions
of Algebraic Geometry and Algebra and Number Theory for the
interest to his work. The author is also grateful to his
colleagues in Algebraic Geometry group at the University of
Liverpool for the creative atmosphere and general support.

%%%%%%%%%%%%%%%%%%%%%%%%%%%%%%%%%%%%%%%%%%%%%%%%%%%%%%%%%%%%%%%%
%%%%%%%%%%%%%%%%%%%%%%%%%%%%%%%%%%%%%%%%%%%%%%%%%%%%%%%%%%%%%%%%
%%%%%%%%%%%%%%%%%%% SECTION 1

\section{Hypersurfaces with a large singular locus}

In this section we consider the problem of estimating the
codimension of the set of polynomials $f$, defining hypersurfaces
with a large singular locus. As a first example we compute the
codimension of the set of polynomials $f$, such that the
hypersurface $\{f=0\}$ has a linear subspace of singular points
(Subsection 1.1). The next by complexity case, when the singular
locus is a hypersurface in a linear subspace, is made in
Subsection 1.2. In Subsection 1.3 we give a precise setting of the
problem in the general case and state the main
result.\vspace{0.3cm}

{\bf 1.1. Hypersurfaces with a linear subspace of singular
points.} Let ${\mathbb P}^N$ be the complex projective space with
homogeneous coordinates $(x_0 : x_1 : \cdots : x_N)$, $N\geqslant
3$ and
$$
{\cal P}_{d,N+1}=H^0({\mathbb P}^N, {\cal O}_{{\mathbb P}^N}(d))
$$
the linear space of homogeneous polynomials of degree $d$ in
$x_*$. For a polynomial $f\in {\cal P}_{d,N+1}\setminus\{0\}$ the
set of singular points of the hypersurface $\{f=0\}$ is denoted by
the symbol $\mathop{\rm Sing}(f)$. Set
$$
{\cal P}^{(i)}_{d,N+1}=\{f\in{\cal P}_{N,d}\,|\,\mathop{\rm
dim}\mathop{\rm Sing}(f)\geqslant i\},
$$
where for the identically zero polynomial $f\equiv 0$ we set
$\mathop{\rm Sing}(0)={\mathbb P}^N$. Obviously, the sets ${\cal
P}^{(i)}_{d,N+1}$ are closed and for $i\geqslant j$ we have ${\cal
P}^{(j)}_{d,N+1}\subset {\cal P}^{(i)}_{d,N+1}$.\vspace{0.1cm}

In order work with the sets ${\cal P}^{(i)}_{d,N+1}$, it is
convenient to represent them as a union of more special subsets
that take into account more information about the set $\mathop{\rm
Sing}(f)$, not only its dimension. For $k\geqslant i$ let
$$
{\cal P}^{(i,k)}_{d,N+1}\subset {\cal P}^{(i)}_{d,N+1}
$$
be the closure of the set ${\cal P}^{(i)}_{d,N+1}$, consisting of
polynomials $f$, such that $\mathop{\rm Sing}(f)$ contains an
irreducible component $C$ of dimension $i\geqslant 1$, the linear
span $\langle C\rangle$ of which is a $k$-plane in ${\mathbb
P}^N$. For instance, ${\cal P}^{(i,i)}_{d,N+1}$ consists of
polynomials $f$, such that $\mathop{\rm Sing}(f)$ contains a
$i$-plane in ${\mathbb P}^N$. The closure in this case is not
needed: the set ${\cal P}^{(i,i)}_{d,N+1}$ allows the following
obvious explicit description.\vspace{0.1cm}

{\bf Proposition 1.1.} {\it The following equality holds:}
$$
\mathop{\rm codim}\left({\cal P}^{(i,i)}_{d,N+1}\subset {\cal
P}_{d,N+1}\right)= {d+i\choose d}+(N-i){d+i-1\choose
d-1}-(i+1)(N-i).
$$

{\bf Proof.} For an $i$-plane $P\subset {\mathbb P}^N$ by the
symbol ${\cal P}^{(i,i)}_{d,N+1}(P)$ we denote the closed set of
polynomials $f$, such that $P\subset \mathop{\rm Sing}(f)$. Fixing
$P$, we may assume that
$$
P=\{x_{i+1}=\dots=x_N=0\},
$$
so that the property $f\in {\cal P}^{(i,i)}_{d,N+1}(P)$ is defined
by the set of identical equalities
$$
f|_P\equiv\left.\frac{\partial f}{\partial
x_{i+1}}\right|_P\equiv\dots\equiv\left.\frac{\partial f}{\partial
x_N}\right|_P\equiv 0.
$$
It is easy to see that these equalities are equivalent to
vanishing the coefficients at the monomials in $x_0,\dots,x_i$ and
the monomials of the form
$$
x_j x^{a_0}_0x^{a_1}_1\dots x^{a_i}_i,
$$
$j=i+1,\dots,N$, which are all distinct. Therefore, the
codimension of the closed set ${\cal P}^{(i,i)}_{d,N+1}(P)$ is
$$
{d+i\choose d}+(N-i){d+i-1\choose d-1}.
$$
Taking into account that for a general polynomial $f\in {\cal
P}^{(i,i)}_{d,N+1}(P)$ the equality $\mathop{\rm Sing}(f)=P$ holds
for some $i$-plane $P$, we obtain the claim of the
proposition.\vspace{0.3cm}

%%%%%%%%%%%%%%%%%%%%%%%%%%%%%%%%%%%%%%%%%%%%%%%%%%%%%%%%%%%%%%%%
%%%%%%%%%%%% subsection 1.2

{\bf 1.2. The singular locus is a hypersurface in a linear
subspace.} Let us consider now the next by complexity example: let
us estimate the codimension of the set ${\cal
P}^{(i,i+1)}_{d,N+1}$. This set is the closure of the set of
polynomials $f$, such that for some $(i+1)$-plane $P\subset
{\mathbb P}^N$ and an irreducible hypersurface $C\subset P$ of
degree $q\geqslant 2$ we have the inclusion $C\subset \mathop{\rm
Sing}(f)$. Fixing the linear subspace $P$, we obtain the closed
subset ${\cal P}^{(i,i+1)}_{d,N+1}(P)\subset {\cal
P}^{(i,i+1)}_{d,N+1}$, so that
$$
{\cal P}^{(i,i+1)}_{d,N+1}=\mathop{\bigcup}\limits_{P\subset
{\mathbb P}^N} {\cal P}^{(i,i+1)}_{d,N+1}(P),
$$
where the union is taken over all $(i+1)$-planes in ${\mathbb
P}^N$. By Bertini's theorem (and the explicit description of the
polynomials $f\in {\cal P}^{(i,i+1)}_{d,N+1}(P)$, given below in
the proof of Proposition 1.2), for a general polynomial $f\in
{\cal P}^{(i,i+1)}_{d,N+1}$ there is a unique $(i+1)$-plane
$P\subset {\mathbb P}^N$, such that $f\in {\cal
P}^{(i,i+1)}_{d,N+1}(P)$, and for that reason the equality
$$
\mathop{\rm codim}\left({\cal P}^{(i,i+1)}_{d,N+1}\subset {\cal
P}_{d,N+1}\right)=\mathop{\rm codim}\left({\cal
P}^{(i,i+1)}_{d,N+1}(P)\subset {\cal
P}_{d,N+1}\right)-(i+2)(N-i-1)
$$
holds.\vspace{0.1cm}

Now let us fix $P$: we may assume that
$$
P=\{x_{i+2}=\dots =x_N=0\}.
$$
It is clear that $C\subset \mathop{\rm Sing}(f|_P)$. If
$f|_P\not\equiv 0$, then $C$ is a multiple component of the
hypersurface $\{f|_P=0\}$. There are at most
$\left[\frac{d}{4}\right]$ such components and they are determined
by the polynomial $f|_P$. However, there is also another option:
$P=\mathop{\rm Sing}(f|_P)$, that is, $f|_P\equiv 0$. In that case
the subvariety of singularities $C\subset P$ is determined by the
polynomial $f$, but not by its restriction $f|_P$. In order to
take both options into account, let us write
$$
{\cal P}^{(i,i+1)}_{d,N+1}(P)={\cal P}^{(i,i+1;i)}_{d,N+1}(P)\cup
{\cal P}^{(i,i+1;i+1)}_{d,N+1}(P),
$$
where ${\cal P}^{(i,i+1;l)}_{d,N+1}(P)$ is the closure of the set
of polynomials $f\in {\cal P}^{(i,i+1)}_{d,N+1}(P)$, such that
$$
\mathop{\rm dim} \mathop{\rm Sing}(f|_P)=l.
$$
The codimension of the set ${\cal P}^{(i,i+1)}_{d,N+1}(P)$ is the
minimum of the codimensions of those two sets. It is obvious from
the explicit formulas for those codimensions that the minimum is
attained at the first set.\vspace{0.1cm}

{\bf Proposition 1.2.} (i) {\it For $d\geqslant 4$, $d\neq 6$ or
for $d=6$, $i\leqslant N-2$ the following equality holds:}
$\mathop{\rm codim}\left({\cal P}^{(i,i+1;i)}_{d,N+1}(P)\subset
{\cal P}_{d,N+1}\right)=$
$$
={d+i+1\choose i+1}-{d+i-3\choose i+1}-{i+3\choose
i+1}+(N-i-1)\left({d+i\choose i+1}-{d+i-2\choose i+1}\right).
$$

(ii) {\it For $d=6$, $i=N-1$ the following equality holds}
$$
\mathop{\rm codim}\left({\cal P}^{(N-1,N;N-1)}_{6,N+1}(P)\subset
{\cal P}_{6,N+1}\right)={N+6\choose 6}-{N+3\choose 3}-1.
$$

(iii) {\it The following equality holds}
 $\mathop{\rm codim}\left({\cal
P}^{(i,i+1;i+1)}_{d,N+1}(P)\subset {\cal P}_{d,N+1}\right)=$
$$
=\mathop{\rm min}\left\{{d+i+1\choose i+1}-{i+3\choose
i+1}+(N-i-1)\left({d+i\choose i+1}-{d+i-2\choose
i+1}\right),\right.
$$
$$
\left.{d+i+1\choose i+1}+(N-i-2)\left({d+i\choose
i+1}-1\right)\right\}.
$$

{\bf Proof.} Let us show the claim (i). For a general polynomial
$f\in {\cal P}^{(i,i+1;i)}_{d,N+1}(P)$ we have: $f|_P\not\equiv
0$, the hypersurface $\{f|_P=0\}\subset P$ has a multiple
component $C$ of degree $q\geqslant 2$, and moreover,
\begin{equation}\label{26.04.2016.1}
\left.\frac{\partial f}{\partial x_j}\right|_C\equiv 0
\end{equation}
for $j=i+2,\dots,N$. Note that the coefficients of the polynomials
$f|_P$ and $\left.\frac{\partial f}{\partial x_j}\right|_P$,
$j=i+2,\dots,N$, corrrespond to {\it distinct} coefficients of the
original polynomial $f$. The requirement that the hypersurface
$\{f|_P=0\}\subset P$ has a double component $C$ of degree
$q\geqslant 2$, gives
$$
E_q={d+i+1\choose i+1}-{d-2q+i+1\choose i+1}-{q+i+1\choose i+1}
$$
independent conditions for the coefficients of the polynomial
$f|_P$.\vspace{0.1cm}

{\bf Lemma 1.1.} {\it For $d\neq 6$ the minimum of the numbers
$E_q$, $q=2,\dots ,[d/2]$, is attained at $q=2$.}\vspace{0.1cm}

{\bf Proof.} We have: $(i+1)!(E_{q+1}-E_q)=$
$$
=[(d-2q+i+1)\dots (d-2q+1)-(d-2q+i-1)\dots (d-2q-1)]-
$$
$$
-[(q+i+2)\dots (q+2)-(q+i+1)\dots (q+1)],
$$
whence after simplifications we get $i! (E_{q+1}-E_q)=$
$$
=(2d-4q+i)(d-2q+i-1)\dots (d-2q+1)-(q+i+1)\dots (q+2).
$$
The first product decreases when $q$ is increasing, the second one
is increasing. It is easy to check that for $d\geqslant 7$ the
inequality $E_3>E_2$ holds. Therefore, for the sequence of
integers $E_q,q=2,\dots,[d/2]$, there are two options:

\begin{itemize}

\item either it is increasing,

\item or it is first increasing $(q=2,\dots,q^*)$, and then
decreasing $(q=q^*,\dots,[d/2])$.

\end{itemize}

In the first case the claim of the lemma is obvious. In the second
case the minimum of the numbers $E_q$ is attained either at $q=2$,
or at $q=[d/2]$, and an easy check shows that the minimum
corresponds precisely to the value $q=2$. Q.E.D. for the
lemma.\vspace{0.1cm}

Now let us fix the polynomial $f|_P\not\equiv 0$. Since the
hypersurface $\{f|_P=0\}$ has finitely many components, we may
assume that the irreducible hypersurface $C$ of degree $q\geqslant
2$ is fixed. Now the requirement (\ref{26.04.2016.1}) imposes on
the coefficients of the polynomial $({\partial f}/{\partial
x_j})|_P$ precisely
\begin{equation}\label{27.04.2016.1}
{d+i\choose i+1}-{d-q+i\choose i+1}
\end{equation}
independent conditions, and it is obvious, that the minimum of the
last expression is attained at $q=2$. This completes the proof of
the claim (i) (an explicit check shows that it is true for $d=6$,
too, although for $d=6$ the claim of Lemma 1.1 is not true:
$E_3<E_2$.) The claim (ii) is shown by explicit simple
computations.\vspace{0.1cm}

Let us show the claim (iii). In that case the hypersurface
$\{f=0\}$ contains the entire subspace $P$. The closed subset of
polynomials $f$, such that $f|_P\equiv 0$, has codimension
${d+i+1\choose i+1}$. Furthermore, either all partial derivatives
$\frac{\partial f}{\partial x_j}, j=i+2.\dots,N$, vanish on $P$,
so that $P\subset \mathop{\rm Sing}(f)$ and this gives an
essentially higher codimension than what is claimed by (iii), and
for that reason this option can be ignored, or $\frac{\partial
f}{\partial x_j}|_P\not\equiv 0$ for some $j\in\{i+2,\dots,N\}$.
Without loss of generality we may assume that $j=i+2$. Then all
polynomials $\partial f/\partial x_j$, $j=i+3,\dots,N$, vanish on
one of the components $C$ of the hypersurface
$\displaystyle\left\{\frac{\partial f}{\partial
x_{i+2}}=0\right\}\subset P$, where $\mathop{\rm deg}C=q\geqslant
2$. This component can be assumed to be fixed and this gives for
each of the $(N-i-2)$ polynomials $\partial f/\partial x_j$,
$j=i+3,\dots,N$, the new independent conditions, the number of
which is given by the formula (\ref{27.04.2016.1}). Taking into
account that the hypersurface $\{\partial f/\partial x_{i+2}=0\}$
is reducible, we finally obtain
$$
{d+i+1\choose i+1}+(N-i-1)\left({d+i\choose i+1}-{d-q+i\choose
i+1}\right)-{i+q+1\choose i+1}
$$
independent conditions for the coefficients of the polynomial $f$.
Using the same method as in the proof of the claim (i), it is easy
to show that the minimum of the last expression for
$q=2,\dots,d-1$ is attained at one of the end values of $q$:
either at $q=2$, or at $q=d-1$. Proof of Proposition 1.2 is
complete. Q.E.D.\vspace{0.3cm}

%%%%%%%%%%%%%%%%%%%%%%%%%%%%%%%%%%%%%%%%%%%%%%%%%%%%%%%%%%%%%%%%
%%%%%%%%%%%% subsection 1.3

{\bf 1.3. Statement of the main result for the general case.} Now
let us consider the general case. By the symbol ${\cal
P}^{(i,k;l)}_{d,N+1}(P)$ for a $k$-plane $P\subset {\mathbb P}^N$
denote the closure of the set of polynomials $f\in {\cal
P}_{d,N+1}$ such that:

\begin{itemize}

\item the set $\mathop{\rm Sing}(f)$ has an irreducible component
$C\subset P$ of dimension $i$, and moreover $\langle C\rangle=P$,

\item the set $\mathop{\rm Sing}(f|_P)$ has an irreducible component
$B$ of dimension $l\geqslant i$, and moreover $C\subset B$.

\end{itemize}

By the symbol ${\cal P}^{(i,k)}_{d,N+1}(P)$ for a $k$-plane
$P\subset {\mathbb P}^N$ denote the closure of the set of
polynomials $f\in {\cal P}_{d,N+1}$ such that the first of the two
conditions stated above is satisfied. Obviously,
$$
{\cal P}^{(i,k)}_{d,N+1}(P)=\mathop{\bigcup}\limits^k_{l=i} {\cal
P}^{(i,k;l)}_{d,N+1}(P).
$$
Everywhere in the sequel the codimension of various closed sets in
the space of polynomials ${\cal P}_{d,N+1}$ is meant to be with
respect to that space, so that, for instance, $\mathop{\rm
codim}{\cal P}^{(i,k)}_{d,N+1}(P)$ is the minimum of the
codimensions $\mathop{\rm codim}{\cal P}^{(i,k;l)}_{d,N+1}(P)$,
where $l=i,\dots,k$, and the following estimate holds:
$$
\mathop{\rm codim}{\cal P}^{(i,k)}_{d,N+1}\geqslant \mathop{\rm
codim}{\cal P}^{(i,k)}_{d,N+1}(P)-(k+1)(N-k).
$$

{\bf Remark 1.1.} It is easy to see that for $N-k<l-i$ the
singular set of a polynomial $f\in {\cal P}^{(i,k;l)}_{d,N+1}(P)$
is of dimension at least $i+1$. Indeed, $\mathop{\rm Sing}(f|_P)$
contains an $l$-dimensional irreducible component $C\subset P$,
where $l>i$. If $k=N$, then $C\subset \mathop{\rm Sing}(f)$. If
$k<N$, then
$$
\left[C\cap\left\{ \left.\frac{\partial f}{\partial
x_{k+1}}\right|_P=\dots =\left.\frac{\partial f}{\partial
x_N}\right|_P=0\right\}\right]\subset \mathop{\rm Sing}(f),
$$
so that $\mathop{\rm dim} \mathop{\rm Sing}(f)\geqslant
l-(N-k)>i$, as we claimed. For that reason everywhere below,
whenever we consider the sets ${\cal P}^{(i,k;l)}_{d,N+1}(P)$, we
assume that the inequality $N+i\geqslant k+l$ holds.\vspace{0.1cm}

In order to give a compact statement of the main result about
these codimensions, we introduce one notation more. For positive
integers $a,b,c$ set
$$
\tau(a,b,c)=\mathop{\rm max} \left\{{a+c\choose c},\,\,
ab+1\right\}.
$$
If we fix $c$, then the first of the two numbers exceeds the
second one for the values of $a$ that are higher than a number of
order $\displaystyle\frac{c}{e}N^{\frac{1}{c}}$, where $e$ is the
base of the natural logarithm. If we fix $a$, the first of the two
numbers exceeds the second one for the values of $c$ that are
higher than a number of order
$\displaystyle\frac{a}{e}N^{\frac{1}{a}}$. The meaning of the
function $\tau$ will be clear in \S 2.\vspace{0.1cm}

{\bf Theorem 1.1.} {\it For $l\leqslant k-2$ the following
estimate holds:}
$$
\mathop{\rm codim}{\cal P}^{(i,k;l)}_{d,N+1}(P)\geqslant
(k-l+1){d+l-2\choose l+1}+(N+i-k-l)\,\tau(d-1,k,i).
$$

{\bf Theorem 1.2.} (i) {\it For $l=k-1$, $d\neq 6$ the following
estimate holds:}
$$
\mathop{\rm codim}{\cal P}^{(i,k;l)}_{d,N+1}(P)\geqslant
\left[{d+k\choose k}-{d-4+k\choose k}-{k+2\choose k}\right]
+(N+i-2k+1)\,\tau(d-1,k,i).
$$

(ii) {\it For $l=k-1, d=6$ the following estimate holds:}
$$
\mathop{\rm codim}{\cal P}^{(i,k;l)}_{d,N+1}(P)\geqslant
\left[{k+6\choose k}-{k+3\choose k}\right]
+(N+i-2k+1)\,\tau(5,k,i).
$$

{\bf Theorem 1.3.} {\it For $l=k$ the following estimate holds:}
$$
\mathop{\rm codim}{\cal P}^{(i,k;l)}_{d,N+1}(P)\geqslant
{d+k\choose k}+(N+i-k-l)\,\tau(d-1,k,i).
$$

{\bf Proof} of these three theorems will be given in \S 2.

%%%%%%%%%%%%%%%%%%%%%%%%%%%%%%%%%%%%%%%%%%%%%%%%%%%%%%%%%%%%%%%%
%%%%%%%%%%%%%%%%%%%%%%%%%%%%%%%%%%%%%%%%%%%%%%%%%%%%%%%%%%%%%%%%
%%%%%%%%%%%%%%%%%%% SECTION 2

\section{Good sequences and linear spans}

In this section we prove Theorems 1.1-1.3. Theorem 1.1 is the
hardest one. First (Subsection 2.1) we describe the strategy of
the proof of this theorem and give the definition of good
sequences and associated subvarieties: this technique amkes it
possible to reconstruct the subvariety $C\subset \mathop{\rm
Sing}(f)$ inside the, generally speaking, larger subvariety
$B\subset \mathop{\rm Sing}(f|_P)$. After that we estimate the
codimension of the subset of polynomials on $P$ with a
non-degenerate subvariety of singular points of dimension $l$
(Subsection 2.2). Finally, in Subsection 2.3 we complete the proof
of Theorem 1.1 by means of well known methods of estimating the
codimension of the set of polynomials, vanishing on a given
non-degenerate subvariety; after that we show Theorems 1.2 and
1.3, which is easy.\vspace{0.3cm}

{\bf 2.1. Plan and start of the proof of Theorem 1.1.} Let us
describe the strategy of the proof of Theorem 1.1. Fix a $k$-plane
$P\subset {\mathbb P}^N$. We will assume that it is the coordinate
plane
$$
P=\{x_{k+1}=\dots =x_N=0\}.
$$
As we noted in Subsection 1.2, the coefficients of the polynomials
$f|_P$ and $\left.\frac{\partial f}{\partial x_j}\right|_P$,
$j=k+1,\dots, N$, correspond to distinct coefficients of the
polynomial $f$. For that reason, considering the general
polynomial $f\in {\cal P}^{(i,k;l)}_{d,N+1}(P)$, one has to solve
three problems:\vspace{0.1cm}

1) estimate the codimension of the closed subset ${\cal
P}^{(l,k)}_{d,k+1}$ in the space ${\cal P}_{d,k+1}$ (since,
obviously, $f|_P\in {\cal P}^{(l,k)}_{d,k+1}$),\vspace{0.1cm}

2) using the $l-i$ polynomials $\left.\frac{\partial f}{\partial
x_j}\right|_P$, where $j\in I\subset\{k+1,\dots,N\}$, so that
$|I|=l-i$, reconstruct the variety of singular points $C\subset
\mathop{\rm Sing}(f)$ as a subvariety of codimension $l-i$ of the
variety of singular points $B\subset \mathop{\rm Sing}(f|_P)$,
which depends on the restriction $f|_P$ only and for that reason
can be assumed to be fixed, if we fix the polynomial $f|_P\in
{\cal P}^{(l,k)}_{d,k+1}$,\vspace{0.1cm}

3) estimate the codimension of the closed set of polynomials $h\in
{\cal P}_{d-1,k+1}$, vanishing on a fixed non-degenerate
subvariety $B\subset P$, and apply this estimate to the
$(N+i-k-l)$ polynomials $\left.\frac{\partial f}{\partial
x_j}\right|_P$, $j\in \{k+1,\dots,N\}$, $j\not\in
I$.\vspace{0.1cm}

The sum of the estimate, obtained at the stage 1), with the
$(N+i-k-l)$-multiple of the estimate, obtained at the stage 3), is
precisely the inequality, claimed by Theorem 1.1.\vspace{0.1cm}

Let us start to realize this programme.\vspace{0.1cm}

First of all, recall the following definition (see \cite[Section
3]{Pukh01} or \cite[Chapter 3]{Pukh13a}).\vspace{0.1cm}

{\bf Definition 2.1.} A sequence of homogeneous polynomials
$g_1,\dots,g_m$ of arbitrary degrees on the projective space
${\mathbb P}^e$, $e\geqslant m+1$, is said to be a {\it good
sequence}, and an irreducible subvariety $W\subset{\mathbb P}^e$
of codimension $m$ is its {\it associated subvariety}, if there
exists a sequence of irreducible subvarieties $W_j\subset{\mathbb
P}^e$, $\mathop{\rm codim}W_j=j$ (in particular, $W_0={\mathbb
P}^e$) such that:

\begin{itemize}

\item $g_{j+1}|_{W_j}\not\equiv 0$ for $j=0,\dots,m+1$,

\item $W_{j+1}$ is an irreducible component of the closed
algebraic set $g_{j+1}|_{W_j}=0$,

\item $W_m=W$.

\end{itemize}

\noindent A good sequence can have more than one associated
subvarieties, but their number is bounded from above by a
constant, depending on the degrees of the polynomials $g_j$ only
(see \cite[Section 3]{Pukh01}).\vspace{0.1cm}

Assuming the polynomial $f|_P$ and the subvariety $B$ to be fixed,
let us construct a good sequence of polynomials on $P={\mathbb
P}^k$ with the subvariety $C$ as one of its associated
subvarieties. This sequence starts with $g_1=f|_P\not\equiv 0$.
Since $B$ is an irreducible $l$-dimensional component of the
closed set $\mathop{\rm Sing}(f|_P)$, for some $(k-l-1)$
polynomials $\left.\frac{\partial f}{\partial x_j}\right|_P$,
$j\in \{0,\dots, k\}$, we obtain a good sequence of polynomials
with $B$ as one of its associated subvarieties. If $l=i$, then
there is nothingmore to construct. Assume that $l\geqslant i+1$.
Then among the polynomials $\left.\frac{\partial f}{\partial
x_j}\right|_P$, $j\in \{k+1,\dots, N\}$, there is one which does
not vanish on $B$ (otherwise, $B\subset \mathop{\rm Sing}(f)$, so
that $f\in {\cal P}^{(l,k;l)}_{d,N+1}(P)$, and this contradicts to
the assumption that $C\subset B$, $C\neq B$ is an irreducible
component of the set $\mathop{\rm Sing}(f)$). We add this
polynomial to already constructed sequence. Continuing in this way
(using at every step the fact that $C$ is an irreducible component
of the set $\mathop{\rm Sing}(f)$), we complete the construction
of a good sequence. Assuming the polynomials of the good sequence
to be fixed, we may assume that the variety $C\subset P$ to be
fixed as well. This solves the problem 2), stated above. Now let
us consider the most difficult problem 1).\vspace{0.3cm}

%%%%%%%%%%%%%%%%%%%%%%%%%%%%%%%%%%%%%%%%%%%%%%%%%%%%%%%%%%%%%%%%
%%%%%%%%%%%% subsection 2.2

{\bf 2.2. Linearly independent singular points.} The problem 1),
stated above, is solved in the following claim.\vspace{0.1cm}

{\bf Proposition 2.1.} {\it The following estimate holds:}
\begin{equation}\label{02.05.2016.1}
\mathop{\rm codim}\left({\cal P}^{(l,k)}_{d,k+1}\subset {\cal
P}_{d,k+1}\right)\geqslant (k-l+1){d+l-2\choose l+1}.
\end{equation}

{\bf Proof.} In order to simplify the notations, we assume that
$k=N$. Let us describe the technique of estimating the codimension
of the closed subset of the space ${\cal P}_{d,N+1}$, consisting
of polynomials with many singular points. The following claim is
true.\vspace{0.1cm}

{\bf Lemma 2.1.} {\it Assume that $d\geqslant 3$. For any set of
$m$ linearly independent points $p_1,\dots,p_m\in{\mathbb P}^N$,
$m\leqslant N+1$, the condition
$$
\{p_1,\dots,p_m\}\subset\mathop{\rm Sing}(g),
$$
$g\in{\cal P}_{d,N+1}$, defines a linear subspace of codimension
$m(N+1$) in ${\cal P}_{d,N=1}$.}\vspace{0.1cm}

{\bf Proof.} We may assume that
$$
p_1=(1:0:0\dots:0),\quad p_2=(0:1:0:\dots:0)
$$
and so on correspond to the first $m$ vectors of the standard
basis of the linear space ${\mathbb C}^{N+1}$. The condition
$p_i\in\mathop{\rm Sing}(g)$ means vanishing of the coefficients
at the monomials $x^d_{i-1},x^{d-1}_{i-1}x_j$, for all $j\neq
i-1$. For $d\geqslant 3$ all these $m(N+1)$ monomials are
distinct. Q.E.D. for the lemma.\vspace{0.1cm}

Now let us consider an arbitrary linear subspace
$\Pi\subset{\mathbb P}^N$ of codimension $r+1$, where $r\geqslant
1$, given by a system of $r+1$ equations
$$
l_0(x)=0,\,\,l_1(x)=0,\dots,l_r(x)=0,
$$
where $l_0,\dots,l_r$ are linearly independent forms. For each
$i=1,\dots,r$ fix an arbitrary set of distinct constants
$\lambda_{i0},\dots,\lambda_{i,d-1}\in{\mathbb C}$; we assume that
$\lambda_{i0}=0$ for all $i=1,\dots,r$. Now for any integer valued
point
$$
\underline{e}=(e_1,\dots,e_r)\in{\mathbb Z}^r_+,\quad e_i\leqslant
d-1,
$$
by the symbol $\Theta(\underline{e})$ we denote the linear
subspace
$$
\{l_i(x)-\lambda_{i,e_i}l_0(x)=0\,|\,i=1,\dots,r\}\subset{\mathbb
P}^N
$$
of codimension $r$. Obviously, $\Theta(\underline{e})\supset\Pi$.
Set
$$
|\underline{e}|=e_1+\dots+e_r\in{\mathbb Z}_+.
$$
For every tuple $\underline{e}\in{\mathbb Z}^r_+$ with
$|\underline{e}|\leqslant d-3$ consider an arbitrary set
$$
S(\underline{e})=\{p_1(\underline{e}),\dots,p_m(\underline{e})\}
\subset\Theta(\underline{e})\backslash\Pi
$$
of $m$ linearly independent points (so that $m\leqslant
N-r+1$).\vspace{0.1cm}

{\bf Proposition 2.2.} {\it The set of conditions
$$
S(\underline{e})\subset\mathop{\rm
Sing(g|_{\Theta(\underline{e})})},
$$
$\underline{e}\in{\mathbb Z}^r_+$, $|\underline{e}|\leqslant d-3$,
defines a linear subspace of codimension
$$
m(N-r+1)|\Delta|
$$
in ${\cal P}_{d,N+1}$, where
$$
\Delta=\{e_1\geqslant 0,\dots,e_r\geqslant
0,e_1+\dots+e_r\leqslant d-3\}\subset{\mathbb R}^r
$$
is an integer valued simplex and $|\Delta|$ is the number of
integral points in that simplex,
$|\Delta|=\sharp(\Delta\cap{\mathbb Z}^r)$.}\vspace{0.1cm}

{\bf Proof.} We may assume that $l_0=x_0$,
$l_1=x_1,\dots,l_r=x_r$. In order to simplify the formulas we will
prove the affine version of the proposition: set
$v_1=x_1/x_0,\dots, v_r=x_r/x_0$ and $u_i=x_{r+i}/x_0$,
$i=1,\dots,N-r$. In the affine space ${\mathbb A}^N\subset{\mathbb
P}^N$, ${\mathbb A}^N={\mathbb P}^N\backslash\{x_0=0\}$ with
coordinates $(u,v)=(u_1,\dots,u_{N-r}$, $v_1,\dots,v_r$) the
affine spaces
$A(\underline{e})=\Theta(\underline{e})\backslash\Pi$ are
contained entirely:
$$
A(\underline{e})=\Theta(\underline{e})\cap{\mathbb A}^N,
$$
so that $S(\underline{e})\subset A(\underline{e})$ for all
$\underline{e}$. Obviously,
$$
A(\underline{e})=\{v_1=\lambda_{1,e_1},\dots,v_r=\lambda_{r,e_r}\}\subset{\mathbb
A}^N
$$
is a $(N-r)$-plane, which is parallel to the coordinate
$(N-r)$-plane $(u_1,\dots,u_{N-r},0,\dots,0)$. Now let us write
the polynomial $g$ in terms of the affine coordinates $(u,v)$ in
the following way:
$$
g(u,v)=\sum_{\underline{e}\in{\mathbb
Z}^r_+,|\underline{e}|\leqslant d}
g_{e_1,\dots,e_r}(u)\prod^r_{i=1}\prod^{e_i-1}_{j=0}(v_i-\lambda_{ij})
$$
(if $e_i=0$, then the corresponding product is meant to be equal
to 1). Here $g_{\underline{e}}(u)=g_{e_1,\dots,e_r}(u)$ is an
affine polynomial in $u_1,\dots,u_{N-r}$ of degree $\mathop{\rm
deg}g_{\underline{e}}\leqslant d-|e|$. When $\lambda_{ij}$ are
fixed, this expression is unique. By Lemma 2.1, the condition
$$
S(\underline{0})=S(0,\dots,0)\subset\mathop{\rm
Sing}(g|_{A(\underline{0})})
$$
defines a linear subspace of codimension $m(N-r+1)$ in the space
of polynomials ${\cal P}_{d,N-r+1}$. However it is easy to see
that
$$
g|_{A(\underline{0})}=g_{0,\dots,0}(u),
$$
since for $\underline{e}\neq\underline{0}$ in the product
$$
\prod^r_{i=1}\prod^{e_i-1}_{j=0}(v_i-\lambda_{ij})
$$
there is at least one factor $(v_i-\lambda_{i0})=v_i$, which
vanishes when we restrict it onto the $(N-r)$-plane
$A(\underline{0})$. Therefore, the condition
$S(\underline{0})\subset \mathop{\rm Sing}(g|_{A(\underline{0})})$
imposes on the coefficients of the polynomial $g_{0,\dots,0}(u)$
precisely $m(N-r+1)$ independent conditions, whereas the
polynomials $g_{\underline{e}}(u)$ for
$\underline{e}\neq\underline{0}$ can be arbitrary.\vspace{0.1cm}

Now let us complete the proof of Proposition 2.2 by induction on
$|\underline{e}|$. More precisely, for every $a\in{\mathbb Z}_+$
set
$$
\Delta_a=\{e_1\geqslant 0,\dots,e_r\geqslant
0,\,\,e_1+\dots+e_r\leqslant a\}\subset{\mathbb R}^r,
$$
so that $\Delta=\Delta_{d-3}$. Let us prove the claim of
Proposition 2.2 in the following form: for every
$a=0,\dots,d-3$\vspace{0.1cm}

\parshape=1
2cm 11cm \noindent $(*)_a$ the set of conditions
$$
S(\underline{e})\subset\mathop{\rm
Sing}(g|_{\Theta(\underline{e})}),
$$
$\underline{e}\in{\mathbb Z}^r_+$, $|\underline{e}|\leqslant a$,
defines a linear subspace of codimension $m(N-r+1)|\Delta_a|$ in
${\cal P}_{d,N+1}$, and, moreover, the linear conditions are
imposed on the coefficients of the polynomials
$g_{\underline{e}}(u)$ for $\underline{e}\in\Delta_a$, whereas for
$\underline{e}\not\in\Delta_a$ the polynomials
$g_{\underline{e}}(u)$ can be arbitrary.\vspace{0.1cm}

The case $a=0$ has already been considered, so that assume that
$a\leqslant d-4$ and the claims $(*)_j$ for $j=0,\dots,a$ have
been shown. Let us show the claim $({*})_{a+1}$. Let
$\underline{e}\in{\mathbb Z}^r_+$ be an arbitrary multi-index,
$|\underline{e}|=a+1$. The restriction onto the affine subspace
$A(\underline{e})$ means the substitution
$v_1=\lambda_{1,e_1},\dots$, $v_r=\lambda_{r,e_r}$. Therefore the
polynomial $g_{\underline{e}(u)}$ comes into the restriction
$g|_{A(\underline{e})}$ with the non-zero coefficient
$$
\alpha_e=\prod^r_{i=1}\prod^{e_i-1}_{j=0}(\lambda_{i,e_i}-\lambda_{ij}).
$$
On the other hand, for $\underline{e}'\neq\underline{e}$,
$|\underline{e}'|\geqslant a+1$ the product
$$
\prod^r_{i=1}\prod^{e'_i-1}_{j=0}(\lambda_{i,e_i}-\lambda_{ij}).
$$
is equal to zero, since for at least one index $i\in\{1,\dots,r\}$
we have $e'_i> e_i$ and so this product contains a factor equal to
zero. Therefore, $g|_{A(\underline{e})}$ is the sum of the
polynomial $\alpha_eg_{\underline{e}}$ and a linear combination of
the polynomials $g_{\underline{e}'}$ with
$|\underline{e}'|\leqslant a$ with constant coefficients. Now,
fixing the polynomials $g_{\underline{e}'}$ with
$|\underline{e}'|\leqslant a$, we obtain that the condition
$$
S(\underline{e})\subset\mathop{\rm Sing}(g|_{A(\underline{e})})
$$
defines an {\it affine} (generally speaking, not linear) subspace
of codimension $m(N-r+1)$ of the space of polynomials
$g_{\underline{e}}(u_1,\dots,u_{N-r})$ of degree at most $d-|e|$,
the corresponding linear subspace of which is given by the
condition
$$
S(\underline{e})\subset\mathop{\rm Sing}g_{\underline{e}}(u).
$$
Moreover, no restrictions are imposed on the coefficients of other
polynomials $g_{\underline{e'}}$ with
$|\underline{e}'|=a+1$.\vspace{0.1cm}

This proves the claim $(*)_a$ for all $a=0,\dots,d-3$. Proof of
Proposition 2.2 is complete.\vspace{0.1cm}

Now let
$$
\Theta=\Theta[l_0,\dots,l_r;\lambda_{i,j},i=1,\dots,r,j=0,\dots,d-1]=
\{\Theta(\underline{e})\,|\,\underline{e}\in\Delta\}
$$
be some set of linear subspaces of codimension $r$ in ${\mathbb
P}^N$, considered in Proposition 2.2. We define the subset
$$
{\cal P}_{d,N+1}(\Theta)\subset{\cal P}_{d,N+1}
$$
by the following condition: for every subspace
$\Theta(\underline{e})$ with $|\underline{e}|\leqslant d-3$ there
is a set
$S(\underline{e})\subset\Theta(\underline{e})\backslash\Pi$,
consisting of $m$ linearly independent points, such that
$S(\underline{e})\subset\mathop{\rm
Sing}(g|_{\Theta(\underline{e})})$.\vspace{0.1cm}

{\bf Proposition 2.3.} {\it The following inequality holds:}
$$
\mathop{\rm codim}({\cal P}_{d,N+1}(\Theta)\subset{\cal
P}_{d,N+1})\geqslant m|\Delta|.
$$

{\bf Proof} is obtained by means of the obvious dimension count:
the subspaces $\Theta(\underline{e})$ are fixed, so that every
point $p_i(\underline{e})$ varies in a $(N-r)$-dimensional family.
Q.E.D. for the proposition.\vspace{0.1cm}

Finally, let us complete the proof of Proposition 2.1. Set $N=k$,
so that the space of polynomials of degree $d$ is ${\cal
P}_{d,k+1}$. For an arbitrary set
$\Theta=\{\Theta(\underline{e})\,|\, \underline{e}\in\Delta\}$ of
linear subspaces of codimension $l$ in $P={\mathbb P}^k$ let
$$
{\cal P}^{(l,k)}_{d,k+1}(P,\Theta)\subset{\cal P}^{(l,k)}_{d,k+1}
$$
be the set of polynomials $h\in{\cal P}^{(l,k)}_{d,k+1}$ such that
the set $\mathop{\rm Sing}(h)$ has an irreducible component $Q$ of
dimension $l$, where $\langle Q\rangle=P$ and the variety $Q$ is
in general position with the subspaces from the set $\Theta$: for
all $\underline{e}\in\Delta$ the set $\Theta(\underline{e})\cap Q$
contains $(k-l+1)$ linearly independent points. Since $\langle
Q\rangle=P$, the subset ${\cal P}^{(l,k)}_{d,k+1}(P,\Theta)$ is a
Zariski open subset of the set ${\cal P}^{(l,k)}_{d,k+1}$, so that
the inequality (\ref{02.05.2016.1}) will be shown, if we prove it
for ${\cal P}^{(l,k)}_{d,k+1}(\Theta)$ instead of ${\cal
P}^{(l,k)}_{d,k+1}$. However, for ${\cal
P}^{(l,k)}_{d,k+1}(\Theta)$ this inequality follows immediately
from Proposition 2.3, since in that case $m=k-l+1$ and
$|\Delta|={d+l-2 \choose l+1}$.\vspace{0.1cm}

Proof of Proposition 2.1 is complete. Q.E.D.\vspace{0.3cm}

%%%%%%%%%%%%%%%%%%%%%%%%%%%%%%%%%%%%%%%%%%%%%%%%%%%%%%%%%%%%%%%%
%%%%%%%%%%%% subsection 2.3

{\bf 2.3. Polynomials, vanishing on a given variety.} By the
symbol ${\cal P}_{d-1,k+1}(B)$ we denote the closed subset of
polynomials $h\in {\cal P}_{d-1,k+1}$ such that $h|_B\equiv 0$ for
a fixed irreducible subvariety $B$. In our case $\mathop{\rm dim}
B=i$ and $\langle B\rangle={\mathbb P}^k$. There are two methods
of estimating the codimension of the set ${\cal
P}_{d-1,k+1}(B)$.\vspace{0.1cm}

The first method was developed in \cite{Pukh98b} (see also
\cite[Chapter 3]{Pukh13a}). Consider a general linear projection
${\mathbb P}^k\dashrightarrow {\mathbb P}^i$, so that $\pi|_B$ is
a regular surjective map. For any non-zero polynomial $g\in {\cal
P}_{d-1,i+1}$ we have $\pi^* g|_B\not\equiv 0$, so that
$$
\mathop{\rm codim}{\cal P}_{d-1,k+1}(B)\geqslant {\cal
P}_{d-1,i+1}={d-1+i\choose i}.
$$

The second method was developed in \cite{Pukh01} (see also
\cite[Chapter 3]{Pukh13a}). Since $\langle B\rangle={\mathbb
P}^k$, a non-zero linear form can not vanish on $B$. Therefore,
the closed subset of decomposable forms
$$
\underbrace{{\cal P}_{1,k+1}\cdot {\cal P}_{1,k+1}\cdot \cdots
\cdot {\cal P}_{1,k+1}}_{d-1}\subset {\cal P}_{d-1,k+1}
$$
intersects with ${\cal P}_{d-1,k+1}(B)$ by zero only, so that
$$
\mathop{\rm codim}{\cal P}_{d-1,k+1}(B)\geqslant (d-1)k+1.
$$
Finally we get:
$$
\mathop{\rm codim}({\cal P}_{d-1,k+1}(B)\subset {\cal
P}_{d-1,k+1})\geqslant \tau(d-1,k,i).
$$
By the arguments of Subsection 2.1 about good sequences and
Proposition 2.1, this completes the proof of Theorem
1.1.\vspace{0.1cm}

Let us show Theorems 1.2 and 1.3. By the arguments given above, we
only need to prove the inequalities
$$
\mathop{\rm codim}({\cal P}^{(k-1,k)}_{d,k+1}\subset {\cal
P}_{d,k+1})\geqslant {d+k\choose k}-{d-4+k\choose k}-{k+2\choose
k}
$$
for $d\neq 6$ and
$$
\mathop{\rm codim}({\cal P}^{(k-1,k)}_{6,k+1}\subset {\cal
P}_{6,k+1})\geqslant {k+6\choose k}-{k+3\choose k}
$$
(both are in fact equalities), which are obtained by repeating the
arguments that were used in the proof of Proposition 1.2 word for
word. The claim of Theorem 1.3 is obvious.

%%%%%%%%%%%%%%%%%%%%%%%%%%%%%%%%%%%%%%%%%%%%%%%%%%%%%%%%%%%%%%%%
%%%%%%%%%%%%%%%%%%%%%%%%%%%%%%%%%%%%%%%%%%%%%%%%%%%%%%%%%%%%%%%%
%%%%%%%%%%%%%%%%%%% SECTION 3

\section{The codimension of the set of non-factorial \\
hypersurfaces}

In this section we prove Theorem 0.1. First we give the list of
possible values of the parameters $i=N-4, k, l$, so that the
complement to the set of polynomials ${\cal P}^{\rm fact}_{d,N+1}$
is contained in the union of the sets ${\cal
P}^{(i,k;l)}_{d,N+1}$. After that, the estimates obtained in \S\S
1-2 are applied in order to estimate the codimension of the set of
non-factorial hypersurfaces (Subsection 3.1); this estimate is
used for proving part (ii) of Theorem 0.1 in Subsection 3.2 and
part (i) of Theorem 0.1 in Subsection 3.3.\vspace{0.3cm}

{\bf 3.1. Hypersurfaces with the singular locus of codimension
three.} Note the following simple\vspace{0.1cm}

{\bf Proposition 3.1.} {\it For $N\geqslant 7$, $d\geqslant 4$ the
following equality holds: }
$$
\tau(d-1,k,N-4)={d+N-5\choose N-4}.
$$

{\bf Proof:} obvious computations.\vspace{0.1cm}

By Grothendieck's theorem,
$$
{\cal P}_{d,N+1}\setminus {\cal P}^{\rm fact}_{d,N+1}\subset {\cal
P}^{(N-4)}_{d,N+1},
$$
so that in order to estimate the codimension of the set of
non-factorial hypersurfaces, we will estimate the codimension of
the set of hypersurfaces with the singular locus of codimension
three. In the notations of \S 1, the set ${\cal
P}^{(N-4)}_{d,N+1}$ is the union of the following eight sets:
$$
{\cal P}^{(N-4,N-4)}_{d,N+1},\quad {\cal
P}^{(N-4,N-3)}_{d,N+1},\quad {\cal
P}^{(N-4,N-2;N-4)}_{d,N+1},\quad {\cal P}^{(N-4,N-2;N-3)}_{d,N+1},
$$
$$
{\cal P}^{(N-4,N-2;N-2)}_{d,N+1},\quad {\cal
P}^{(N-4,N-1;N-4)}_{d,N+1},\quad {\cal P}^{(N-4,N-1;N-3)}_{d,N+1},
\quad {\cal P}^{(N-4,N)}_{d,N+1}.
$$
We set respectively:
$$
\alpha_1={d+N-4\choose N-4}+4{d+N-5\choose N-4}-4(N-3),
$$
$\alpha_2=\mathop{\rm min} \{\alpha_{2a},\alpha_{2b}\}$, where
$$
\alpha_{2a}={d+N-3\choose N-3}-{d+N-7\choose N-3}+3\left[
{d+N-4\choose N-3}-{d+N-6\choose N-3}\right]-\frac{(N+5)(N-2)}{2},
$$
$$
\alpha_{2b}={d+N-3\choose N-3}+2{d+N-4\choose N-3}-3(N-2).
$$
Furthermore,
$$
\alpha_3=3{d+N-6\choose N-3}+2{d+N-5\choose N-4}-2(N-1),
$$
$$
\alpha_4={d+N-2\choose N-2}-{d+N-6\choose N-2}+{d+N-5\choose
N-4}-\frac{(N+4)(N-1)}{2}
$$
and
$$
\alpha'_4={N+4\choose 6}-{N+1\choose 3}+{N-1\choose 3}-2(N-1)
$$
depends on the dimension $N$ only. Finally,
$$
\alpha_5={d+N-2\choose N-2}-2(N-1),
$$
$$
\alpha_6=4{d+N-6\choose N-3}+{d+N-5\choose N-4}-N,
$$
$$
\alpha_7=3{d+N-5\choose N-2}-N
$$
and
$$
\alpha_8=5{d+N-6\choose N-3}.
$$
Now Propositions 1.1, 1.2 and Theorems 1.1-1.3, taking into
account Proposition 3.1, immediately imply\vspace{0.1cm}

{\bf Theorem 3.1.} {\it For $N\geqslant 7$ the following
inequality holds:}
$$
\mathop{\rm codim}\left({\cal P}_{d,N+1}\setminus {\cal P}^{\rm
fact}_{d,N+1}\right)\geqslant \mathop{\rm min} \{\alpha_i\,|\,
i=1,\dots, 8\}.
$$

{\bf Remark 3.1.} For $d=6$ in this inequality one should replace
$\alpha_4$ by $\alpha'_4$, however, the minimum of the right hand
side is attained at $\alpha_8$ all the same, so that the claim of
Theorem 3.1 remains correct in this case as well.\vspace{0.3cm}

%%%%%%%%%%%%%%%%%%%%%%%%%%%%%%%%%%%%%%%%%%%%%%%%%%%%%%%%%%%%%%%%
%%%%%%%%%%%% subsection 3.2

{\bf 3.2. Hypersurfaces of general type.} In oredr to prove the
claim (ii) of Theorem 0.1, one needs to check that
$\alpha_i\geqslant \alpha_1$ for $i\geqslant 2$, if $d\geqslant
2N$. This check is elementary and we do not perform it here,
giving only one example: setting $d=2N+a$, write
$$
\alpha_8-\alpha_1=4(N-3)+\frac{(3N+a-6)!}{(N-3)!(2N+a)!}\times
$$
$$
\times[10N(2N+a-1)(2N+a-2)-(N-3)(3N+a-5)(11N+5a-4)].
$$
It is easy to see that for $a\geqslant 0$ the expression in the
square brackets is positive, which implies that
$\alpha_8>\alpha_1$. (In fact, the difference $\alpha_8-\alpha_1$
is quite large, but we do not need that.) The remaining
inequalities $\alpha_i>\alpha_1$ for $i\neq 2$ are shown in a
similar way. The reason why $\alpha_1$ realizes the minimum of
$\alpha_i$, $i=1,\dots,8$, is the polynomiality of the functions
$\alpha_i$ in $d$ when the dimension $N$ is fixed: $\alpha_1$ is a
polynomial of the least degree $N-4$. The polynomial $\alpha_{2a}$
also has the degree $N-4$, but its senior coefficient is much
higher. The claim (ii) of Theorem 0.1 is shown.\vspace{0.1cm}

We see that for $d\geqslant 2N$ the irreducible component of the
maximal dimension of the closed set ${\cal P}^{(N-4)}_{d,N+1}$ is
${\cal P}^{(N-4,N-4)}_{d,N+1}$, that is, the set of polynomials
$f\in {\cal P}_{d,N+1}$ such that the hypersurface $\{f=0\}$ has a
$(N-4)$-lane of singular points. \vspace{0.3cm}

%%%%%%%%%%%%%%%%%%%%%%%%%%%%%%%%%%%%%%%%%%%%%%%%%%%%%%%%%%%%%%%%
%%%%%%%%%%%% subsection 3.3

{\bf 3.3. Fano hypersurfaces.} Let us prove the claim (i) of
Theorem 0.1. Again an elementary (but tiresome) check shows that
for $i=1,\dots,6$ the inequality
$$
\alpha_i\geqslant\mathop{\rm min}\{\alpha_7,\alpha_8\}
$$
holds (for $d=6$ with $\alpha_4$ replaced by $\alpha'_4$), which
implies the claim (i). We do not give the tiresome computations
here, except for one example:
$$
\alpha_6-\alpha_8=-N+\frac{(d+N-6)!}{(N-3)!(d-1)!}[-d^2+d(N-1)+
(N^2-9N+18)].
$$
It is easy to see that for $d=4,\dots ,N$ the difference
$\alpha_6-\alpha_8$ is positive. In a similar way the other
inequalities $\alpha_i>\alpha_8$ for $i=1,\dots,5$ are checked.
Proof of the claim (i) of Theorem 0.1 is complete.\vspace{0.1cm}

{\bf Remark 3.2.} Elementary computations, which we do not give
here, show that for $d=4,\dots, d_*(N)$ the inequality
$\alpha_7\leqslant \alpha_8$ holds, and for $d=d_*(N)+1,\dots, N$
the opposite inequality $\alpha_7>\alpha_8$ holds. Here
$d_*(N)\sim \frac23 N$. More precisely, if $N=3m+e$,
$e\in\{0,1,2\}$, then $d_*(N)=2m+e+1$.

%%%%%%%%%%%%%%%%%%%%%%%%%%%%%%%%%%%%%%%%%%%%%%%%%%%%%%%%%%%%%%%%
%%%%%%%%%%%%%%%%%%%%%%%%%%%%%%%%%%%%%%%%%%%%%%%%%%%%%%%%%%%%%%%%
%%%%%%%%%%%%%%%%%%% BIBLIOGRAPHY

\begin{flushleft}
Department of Mathematical Sciences,\\
The University of Liverpool
\end{flushleft}

\noindent{\it pukh@liv.ac.uk}

\end{document}